\documentclass[11pt,a4paper,leqno]{amsart}
\usepackage{amssymb,graphicx,caption,subcaption}
\makeatletter\renewcommand{\@secnumfont}{\relax}\makeatother
\theoremstyle{plain}

\newtheorem{theorem}{Theorem}[section]

\newtheorem{proposition}[theorem]{Proposition}

\theoremstyle{definition}

\theoremstyle{remark}

\newcommand{\proofof}[1]{\end{#1}\begin{proof}} 
\numberwithin{equation}{section}

\begin{document}
\title{On the Atiyah problem on hyperbolic configurations of four points}
\author{Joseph Malkoun}
\address{Department of Mathematics and Statistics\\
Notre Dame University-Louaize, Zouk Mikael\\
Lebanon}
\email{joseph.malkoun@ndu.edu.lb}
\date{\today}
\begin{abstract} Given a configuration $\mathbf{x}$ of $n$ distinct points in hyperbolic $3$-space $H^3$, M.F. Atiyah associated $n$ polynomials 
$p_1,\ldots,p_n$ of a variable $t \in \mathbb{C}P^1$, of degree $n-1$, and conjectured that they are linearly independent over $\mathbb{C}$, 
no matter which configuration $\mathbf{x}$ one starts with. We prove this conjecture for $n=4$ in two cases: 
in case the $4$ points are non-coplanar, and in case one of the points lies in the hyperbolic convex hull of the other three.      
\end{abstract}
\maketitle

\section{Introduction} \label{introduction}

Motivated by physics, particularly by paper \cite{BR1997} by Michael Berry and Jonathan Robbins, 
Michael Atiyah in \cite{Atiyah2000} and \cite{Atiyah2001}, and later on with Paul Sutcliffe in \cite{Atiyah-Sutcliffe2002}, 
wrote a series of papers describing an elementary problem in 
geometry relating configurations of points and polynomials, and made several conjectures. We now describe the hyperbolic version 
of this problem, and the linear independence conjecture in this case. 
We make use of the Poincar\'{e} open $3$-ball model of hyperbolic $3$-space $H^3$. Given 
a configuration $\mathbf{x}=(x_1, \ldots, x_n)$ of $n$ distinct points in hyperbolic $3$-space $H^3$, one 
first forms the points $t_{ij} \in S^2_{\infty}$, 
for $1 \leq i,j \leq n$, $i \neq j$, where $t_{ij}$ is the point on the sphere at infinity obtained as the limit point of 
the hyperbolic ray based at $x_i$ and passing through $x_j$. Identifying $S^2_{\infty}$ with the Riemann sphere $\mathbb{C}P^1$, 
one then defines for each $i$, $1 \leq i \leq n$, the polynomial $p_i$ depending on a variable $t \in \mathbb{C}P^1$, of 
degree $n-1$, with roots $t_{ij}$, for $j \neq i$, each defined up to a non-zero scalar factor. M.F. Atiyah conjectured that 
no matter which configuration $\mathbf{x}$ one starts with, the $n$ polynomials $p_1,\ldots,p_n$ are linearly independent over 
$\mathbb{C}$ (in this paper, we will refer to it as \emph{Atiyah's hyperbolic conjecture}, or simply as AHC). 

\begin{theorem}\label{main-thm} AHC is true for non-coplanar hyperbolic configurations of four points, as well as coplanar hyperbolic configurations of four points with one of the points in the hyperbolic convex hull of the other three. \end{theorem}

We outline the basic idea of our proof here. Details will follow in later sections. It is a proof by contradiction. 
Assume that there is a hyperbolic configuration of $4$ distinct points which belongs to one of the two cases in the theorem, 
and for which the $4$ polynomials $p_1,\ldots,p_4$ are linearly dependent. This implies that there is a linear relation between 
the elementary symmetric polynomials $s_i(t_1,t_2,t_3)$, $0 \leq i \leq 3$ (with $s_0(t_1,t_2,t_3) =1$), which holds when $(t_1,t_2,t_3)$ 
is replaced by any of $(t_{12},t_{13},t_{14})$, $(t_{21},t_{23},t_{24})$, $(t_{31},t_{32},t_{34})$ or $(t_{41},t_{42},t_{43})$. We then make use of 
the orientation-preserving isometry group $PSL(2,\mathbb{C})$ of $H^3$ to simplify this linear relation as much as possible. In the end, 
we show that such a relation cannot occur using the geometry of the problem (namely the fact that the roots are not arbitrary, but rather 
come from a hyperbolic configuration). The last part uses also the Gauss-Lucas theorem, as well as the following result
\begin{proposition} \label{geometric-mean}
Let $D$ be either a closed disk in the complex plane, or a closed half-plane. Prove that for all positive integers $n$, 
and for all complex numbers $z_1,z_2,\ldots,z_n \in D$ there exists a $z \in D$ such that $z^n = z_1 \cdot z_2\cdots z_n$. 
\end{proposition}
This proposition is well known (it was even used as a mathematical olympiad problem in Romania in 2004, in the case where 
$D$ is a closed disk) and 
is an easy consequence of the Grace-Walsh-Szeg\H{o} coincidence theorem.

For the reader's convenience, We state the Gauss-Lucas theorem here, followed by the 
Grace-Walsh-Szeg\H{o} theorem, both without proofs.
\begin{theorem}[Gauss-Lucas] If $P$ is a non-constant polynomial with complex coefficients, all zeros of $P'$ belong to the convex hull of the 
set of zeros of $P'$.
\end{theorem}
Here the notion of convexity is the usual one, in the sense that it uses straight line segments in the 
complex plane.

By a closed circular domain in the complex plane, we mean either a point, or the whole complex plane, or a closed disk 
(meaning the union of a circle and the region ``inside'' the circle), or the complement of an open disk (meaning a circle 
and the region ``outside'' the circle).
\begin{theorem}[Grace-Walsh-Szeg\H{o}] \label{GWS} If $f(z_1,\ldots,z_n)$ is a complex polynomial on $\mathbb{C}^n$ with coordinates 
$z_1, \ldots, z_n$, which is affine in each variable $z_i$ separately, and is symmetric under all permutations of $z_1, \ldots, z_n$, 
and if $\zeta_1,\ldots,\zeta_n$ are $n$ points (distinct or not) in a closed circular domain $D \subseteq \mathbb{C}$, then there is a 
point $\zeta \in D$ such that $f(\zeta_1,\ldots,\zeta_n) = f(\zeta,\ldots,\zeta)$.
\end{theorem}
We remark that applying this theorem to $f(z_1,\ldots,z_n) = z_1 \cdots z_n$ yields proposition \ref{geometric-mean}.

\section{Proof of Theorem \ref{main-thm}}

Assume there exists a configuration of four points $\mathbf{x} = (x_1,\ldots,x_4)$ 
whose associated polynomials $p_1,\ldots,p_4$ are linearly dependent over $\mathbb{C}$. For the moment, 
let us not impose further restrictions on $\mathbf{x}$. Consider 
the $4\times 4$ matrix whose $j$th column is the vector of coefficients of the polynomial $p_j$. By our assumption, 
it is singular, which implies that its rows are linearly dependent over $\mathbb{C}$. Thus there exist complex numbers 
$c_0,\ldots,c_3$ such that
\begin{equation} c_0\cdot s_3(t_1,t_2,t_3) + c_1\cdot  s_2(t_1,t_2,t_3) + c_2 \cdot s_1(t_1,t_2,t_3) + c_3 = 0 \label{relation} \end{equation}
holds simultaneously for $(t_1,t_2,t_3)$ being any of $(t_{12},t_{13},t_{14})$, $(t_{21},t_{23},t_{24})$, $(t_{31},t_{32},t_{34})$ 
or $(t_{41},t_{42},t_{43})$. We will refer to the latter as the four triplets of roots. We remark that, 
for a non-coplanar configuration, the four hyperbolic faces of the hyperbolic 
tetrahedron with vertices $x_1,\ldots,x_4$ define four corresponding distinct circles on $S^2_{\infty}$, by taking limit points 
(the circle ``at infinity'' for that hyperbolic plane). The set of all intersection points of these four circles are 
precisely the twelve roots $t_{ij}$, $i \neq j$ ($1 \leq i,j \leq 4$). We now make use of the complex plane picture to describe 
the Riemann sphere, the latter being the complex plane and a point $\infty$. Each of these circles at infinity must contain exactly
three triplets of roots in the closure of its interior, and one triplet of roots in the closure of its exterior, with each of 
these triplets having exactly two roots on the circle itself, or 
vice versa, meaning the closure of its exterior contains exactly three triplets of roots, with each of these triplets having 
exactly two roots on the circle itself.

\begin{figure}[h!]
\centering
\begin{subfigure}[t]{0.48\textwidth}
  \includegraphics[width=\textwidth]{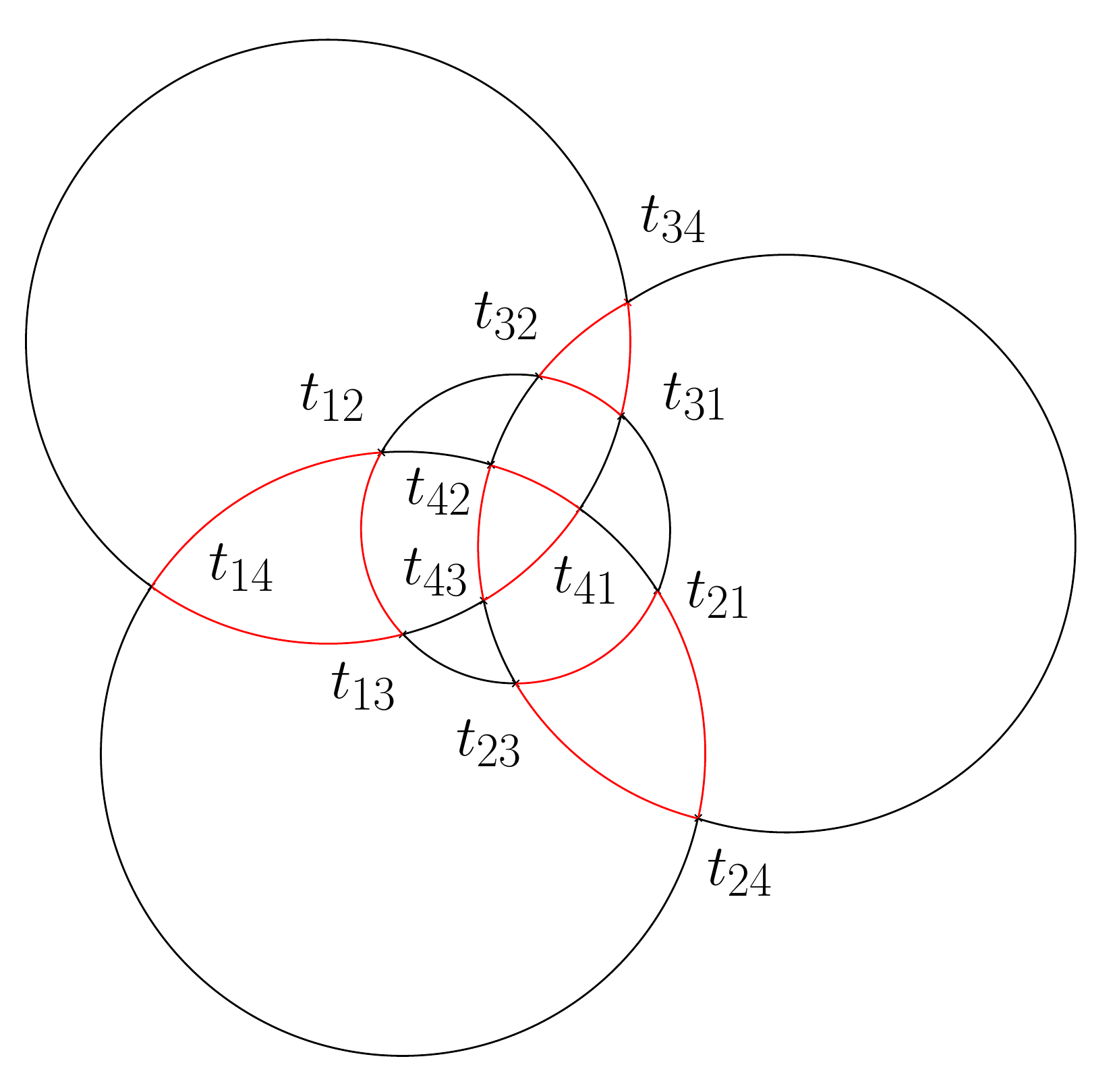}
  \caption*{Type A$1$}
  \label{fig_A1}
\end{subfigure}
\begin{subfigure}[t]{0.48\textwidth}
  \includegraphics[width=\textwidth]{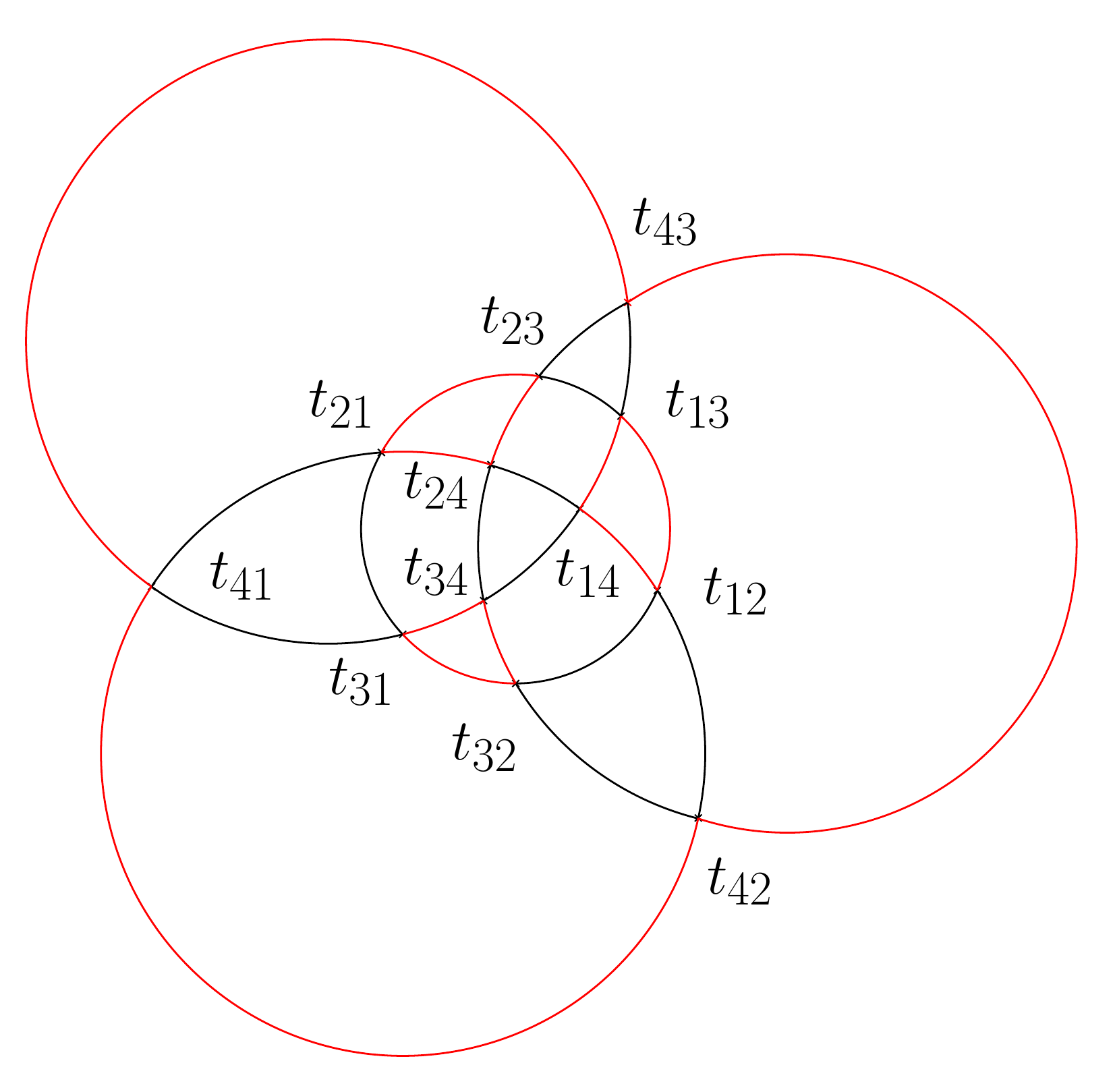}
  \caption*{Type A$2$}
  \label{fig_A2}
\end{subfigure}
\begin{subfigure}[b]{0.48\textwidth}
  \vspace{2 cm}
  \includegraphics[width=\textwidth]{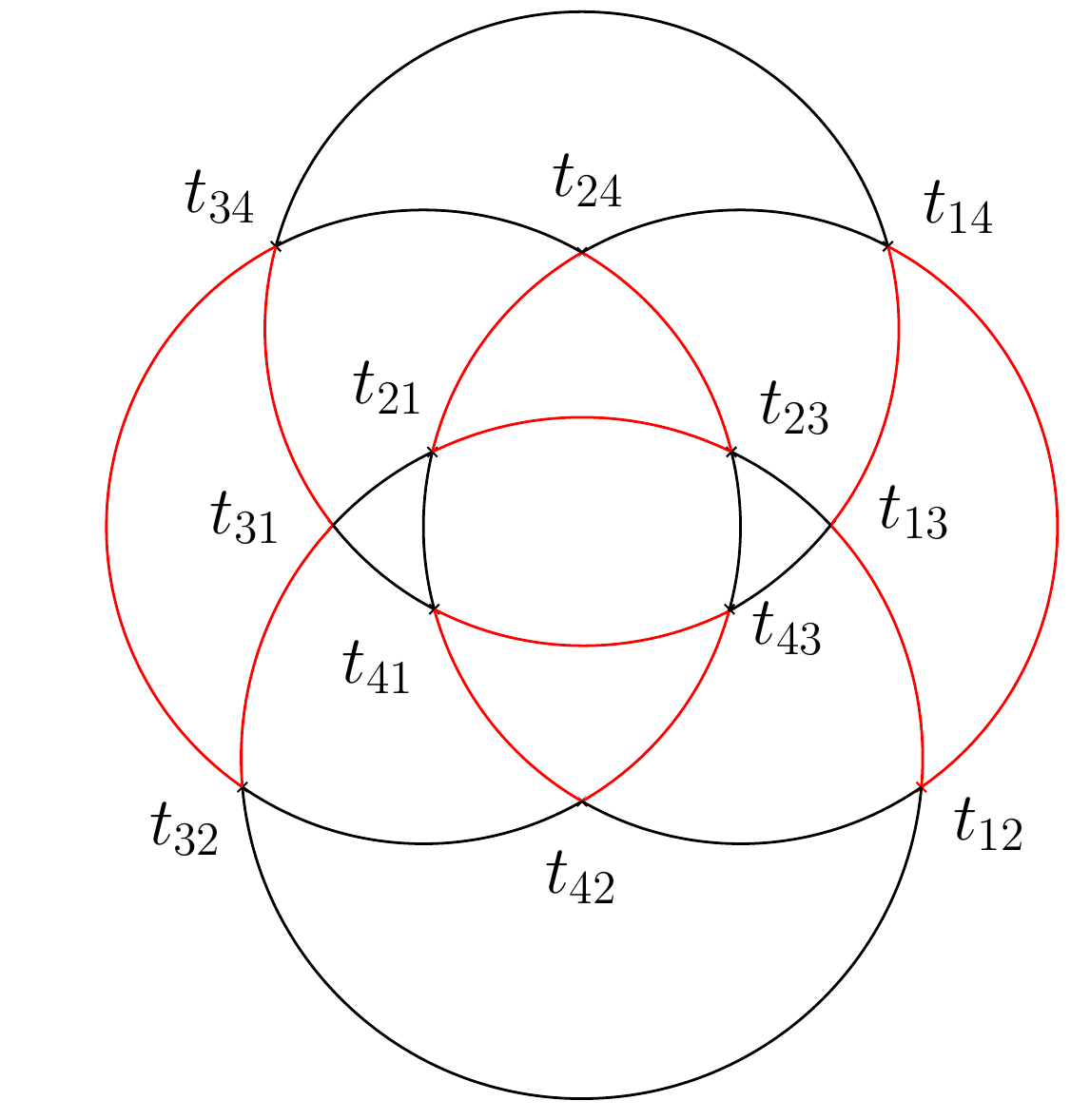}
  \caption*{Type B$1$}
  \label{fig_B1}
\end{subfigure}
\begin{subfigure}[b]{0.48\textwidth}
  \vspace{2 cm}
  \includegraphics[width=\textwidth]{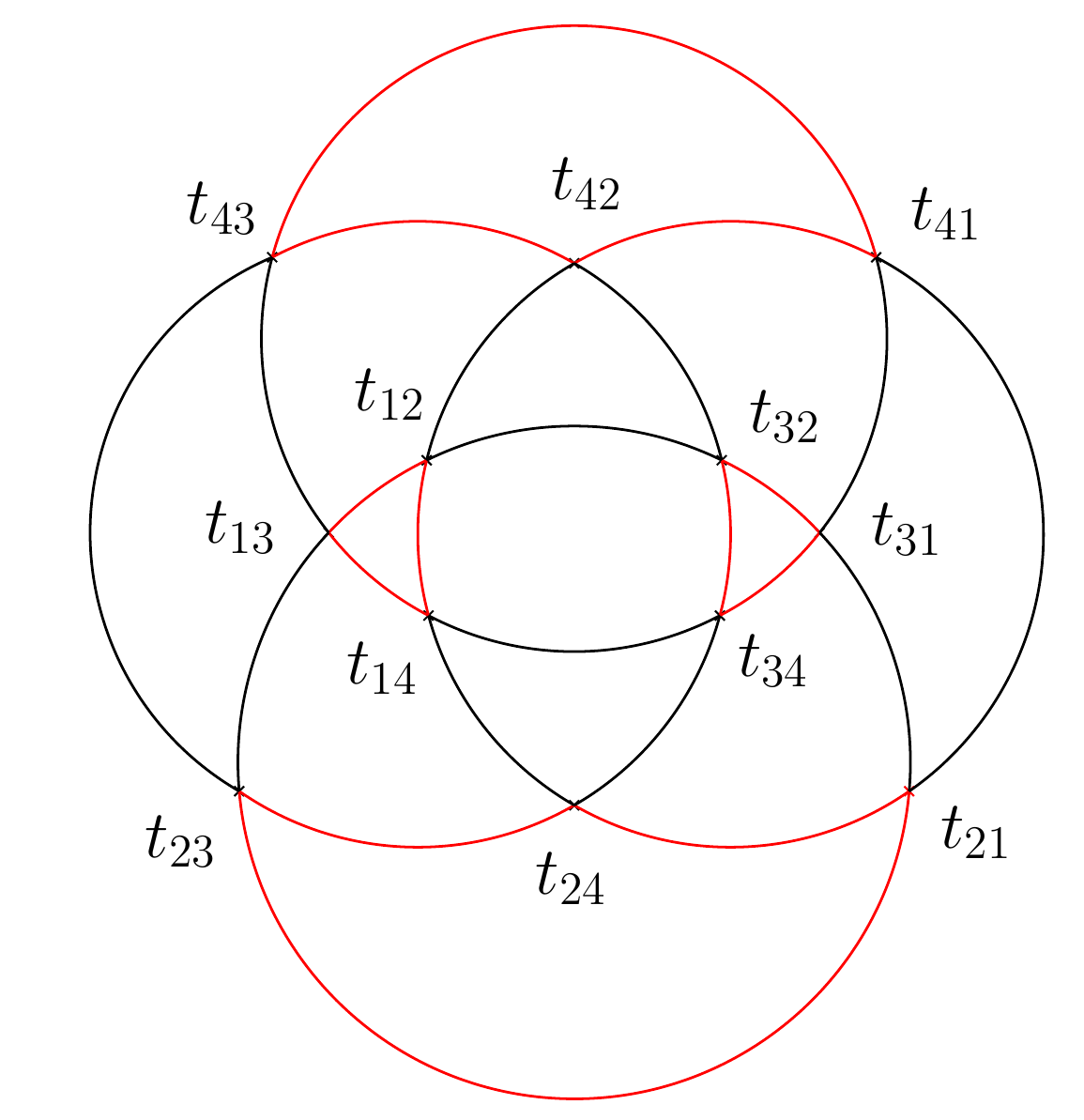}
  \caption*{Type B$2$}
  \label{fig_B2}
\end{subfigure}
\caption{non-coplanar types of configurations of four points}
\label{fig_non-coplanar}
\end{figure}

\begin{figure}[h!]
\centering
\includegraphics[width=0.5\textwidth]{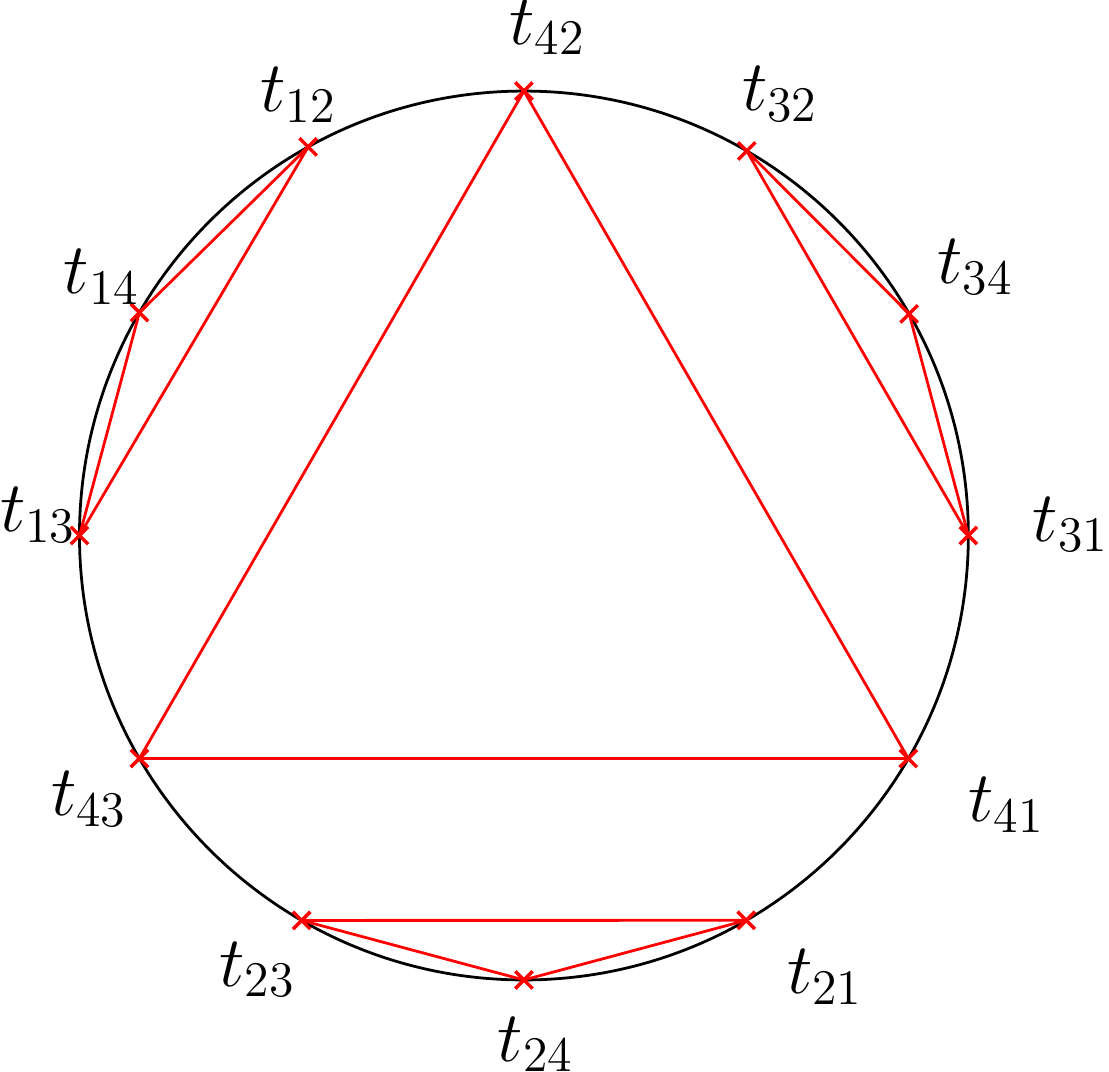}
  \caption{Coplanar case, for a configuration of four points with one of the points in the hyperbolic convex hull of the other three}
  \label{fig_coplanar_special}
\end{figure}

With these constraints in mind, we get only four possible types in the non-coplanar case, described in figure \ref{fig_non-coplanar}, 
up to relabeling the four points $x_1,\ldots,x_4$. In these figures, the triplets of roots, along with the circular arcs connecting the roots, 
are shown in red. On the other hand, figure \ref{fig_coplanar_special} corresponds to the roots of a coplanar configuration of four points, 
with one of the points lying in the hyperbolic convex hull of the other three.

Making use of homogeneous coordinates for $\mathbb{C}P^1$,
\[ w_i = \left( \begin{array}{c} u_i \\
                                 v_i \end{array} \right) \]
for $1 \leq i \leq 3$, such that $t_i = v_i/u_i$ is the corresponding affine coordinate on $\{ u_i \neq 0 \}$, 
we can think of the left-hand side of equation (\ref{relation}) as a symmetric complex trilinear form, which we denote by 
$G(w_1,w_2,w_3)$. With new homogeneous coordinates $w = (u,v)^T$ on $\mathbb{C}P^1$, we define
\[ g(w) = G(w,w,w) \]
Then $G$ can be recovered from the homogeneous cubic polynomial $g$ by the process of complete polarization. Moreover, 
the polarization map $g \mapsto G$ is $PSL(2,\mathbb{C})$-equivariant. On the other hand, the homogeneous cubic polynomial $g$ 
has $3$ roots in $\mathbb{C}P^1$, and we have the following $3$ cases:
\begin{itemize}
\item[a)] $g$ has three distinct roots. Using the $PSL(2,\mathbb{C})$ action, we can without loss of generality assume 
that $g(w) = u v^2 - u^3$, which corresponds to the relation
\begin{equation} \frac{t_1 t_2 + t_2 t_3 + t_3 t_1}{3} = 1 \end{equation}
\item[b)] $g$ has one double root and one single root. Using the $PSL(2,\mathbb{C})$ action, we can assume that 
$g(w) = u^2 v - u^3$, which corresponds to the relation
\begin{equation} \frac{t_1 + t_2 + t_3}{3} = 1 \end{equation}
\item[c)] $g$ has a triple root. Using the $PSL(2,\mathbb{C})$ action, we can assume that 
$g(w) = v^3$, which corresponds to the relation
\begin{equation} t_1 t_2 t_3 = 0 \end{equation}
\end{itemize}

Scenario c) quickly leads to a contradiction. First, we have $t_{12} t_{13} t_{14} = 0$, so 
without loss of generality, assume that $t_{12} = 0$. But we also have $t_{21} t_{23} t_{24} = 0$, 
and clearly $t_{21} \neq t_{12} = 0$, so without loss of generality, assume that $t_{23} = 0$. Thus 
$x_1,x_2$ and $x_3$ all lie on the ray going from $x_1$ towards $0 \in S^2_{\infty}$, and in that order. We have 
that $t_{31} t_{32} t_{34} = 0$, but $t_{31}$ and $t_{32}$ cannot be $0$, so that $t_{34} = 0$. We then have that 
$x_1,x_2,x_3$ and $x_4$ all lie on the ray going from $x_1$ towards $0 \in S^2_{\infty}$, so that, in particular, 
$t_{41}, t_{42}$ and $t_{43}$ are all different from $0$, and the last relation $t_{41} t_{42} t_{43} = 0$ therefore 
leads to a contradiction.  

In scenario b), since $(t_1+t_2+t_3)/3 = 1$, and since the arithmetic mean of $t_1$, $t_2$ and $t_3$ must 
lie in the convex hull of $t_1$, $t_2$ and $t_3$ in the complex plane, it then follows that the convex hulls 
of any two different triplets of roots must intersect non-trivially (since each of them must contain $1$). 
This leads to a contradiction in the two cases 
we are considering, since there are at least three different triplets of roots whose convex hulls are pairwise disjoint 
(and in fact the existence of just two such triplets suffices to give a contradiction).

It remains to consider scenario a). Let 
\[ P(t) = (t-t_1) (t-t_2) (t-t_3) \]
Then
\begin{align*} P'(t) &= 3 t^2 -2 (t_1 + t_2 + t_3) t + t_1 t_2 + t_2 t_3 + t_3 t_1 \\
                     &= 3\left(t^2 - \frac{2}{3}(t_1+t_2+t_3)t+1 \right)
\end{align*}
Let $r_1$ and $r_2$ be the two roots of $P'$. By the Gauss-Lucas theorem, they must lie 
in the convex hull of $t_1$, $t_2$ and $t_3$. We also have $r_1r_2 = 1$. The line segment joining 
$r_1$ with $r_2$ must therefore contain a point on the real line. We deduce that there must be a line 
in the complex plane (the real line) which intersects non-trivially each convex hull of a triplet of roots. 
This is impossible for non-coplanar configurations of types A$1$ and A$2$, as well as coplanar configurations of four points with 
one of the points lying in the hyperbolic convex hull of the other three (figure \ref{fig_coplanar_special}).

There is also another consequence of the existence of $r_1$ and $r_2$ with $r_1 r_2 = 1$ in 
the convex hull of $t_1$, $t_2$ and $t_3$. By theorem \ref{GWS} (or essentially, by proposition \ref{geometric-mean}), 
any circular domain containing $t_1$, $t_2$ and $t_3$ must contain either $1$ or $-1$. But for non-coplanar configurations 
of types B$1$ and B$2$, we can find three disjoint circular domains, each containing a different triplet of roots. This then 
leads to a contradiction, and the theorem is proved.

\def\Dbar{\leavevmode\lower.6ex\hbox to 0pt{\hskip-.23ex \accent"16\hss}D}
\providecommand{\bysame}{\leavevmode\hbox to3em{\hrulefill}\thinspace}
\providecommand{\MR}{\relax\ifhmode\unskip\space\fi MR }
\providecommand{\MRhref}[2]{%
  \href{http://www.ams.org/mathscinet-getitem?mr=#1}{#2}
}
\providecommand{\href}[2]{#2}

\end{document}